\documentstyle[11pt]{article}
\begin{document}

\title{Compatible Poisson brackets \\
of hydrodynamic type}
\author{{\Large Ferapontov E.V. } \\
    Department of Mathematical Sciences \\
    Loughborough University \\
    Loughborough, Leicestershire LE11 3TU \\
    United Kingdom \\
    e-mail: {\tt E.V.Ferapontov@lboro.ac.uk} \\
    and \\
    Centre for Nonlinear Studies \\
    Landau Institute of Theoretical Physics\\
    Academy of Science of Russia, Kosygina 2\\
    117940 Moscow, GSP-1, Russia\\
    e-mail: {\tt fer@landau.ac.ru}
}
\date{}
\maketitle

\newtheorem{theorem}{Theorem}

\pagestyle{plain}

\maketitle

\begin{abstract}

\bigskip

\noindent Some   general  properties of
compatible Poisson brackets of hydrodynamic type are discussed, in
particular:

\medskip

\noindent --- an invariant differential-geometric criterion of the
compatibility
based on the Nijenhuis tensor;

\medskip

\noindent --- the Lax pair with a spectral parameter governing compatible
Poisson 
brackets in the diagonalizable case;

\medskip

\noindent --- the connection of this problem with the class of surfaces
in  Euclidean space which possess nontrivial deformations preserving the
Weingarten operator.

\bigskip

\end{abstract}

\newpage
\section{Introduction}

In 1983 Dubrovin and Novikov \cite{DN83} introduced the Poisson brackets of
hydrodynamic type
\begin{equation}
\{F, G\}=\int \frac{\delta F}{\delta u^i} A^{ij} \frac{\delta G}{\delta u^j}
\ dx
\label{Poisson}
\end{equation}
defined by the Hamiltonian operators $A^{ij}$ of the form
\begin{equation}
A^{ij}=g^{ij}\frac{d}{dx} +b^{ij}_ku^k_x, ~~~~ b^{ij}_k=-g^{is}\Gamma
^j_{sk}.
\label{H1}
\end{equation}
They proved that in the nondegenerate case $(det ~g^{ij}\ne 0)$ the bracket
(\ref{Poisson}), (\ref{H1}) is
skew-symmetric and satisfies the Jacobi identities if and only if the metric
$g^{ij}$ (with upper indices) is flat, and $\Gamma^j_{sk}$
are the Christoffel symbols of the corresponding Levi-Civita connection.

Let us assume that there is a  second Poisson bracket of hydrodynamic type
defined on the same phase
space by the Hamiltonian operator
\begin{equation}
\tilde A^{ij}=\tilde g^{ij}\frac{d}{dx}+\tilde b^{ij}_ku^k_x, ~~~~
\tilde b^{ij}_k=-\tilde g^{is}\tilde \Gamma ^j_{sk},
\label{H2}
\end{equation}
corresponding to a flat metric $\tilde g^{ij}$. Two Poisson brackets
(Hamiltonian operators) are called compatible, if their linear combinations
$\tilde A^{ij}+\lambda A^{ij}$ are Hamiltonian as well.
This requirement implies, in particular, that the metric
$\tilde g^{ij}+\lambda g^{ij}$ must be flat for any  $\lambda$
(plus certain additional restrictions). The necessary and sufficient
conditions of the compatibility
were first formulated by Dubrovin \cite{D93}, \cite{D94} (see
\cite{Mokhov99}, \cite{Mokhov00} for further discussion).
In sect. 2 we reformulate these conditions in terms
of the operator 
$
r^i_j=\tilde g^{is} g_{sj}
$ 
(Theorem 1) which, in particular, imply the vanishing of the Nijenhuis
tensor of the operator $r^i_j$:
$$
N^i_{jk}=r^s_j\partial_sr^i_k  -  r^s_k\partial_sr^i_j  -
r^i_s(\partial_j r^s_k - \partial_kr^s_j)=0
$$ 
(see \cite{Fer90}, \cite{Mokhov00}).

Examples of compatible Hamiltonian pairs naturally arise in the theory of
Hamiltonian systems of hydrodynamic type --- see e.g.
\cite{Arik}, \cite{Gumral90}, \cite{Gumral94}, \cite{Nutku}, \cite{Olver},
 \cite{Pavlov94.1}.
Dubrovin developed a deep theory for a particular class of compatible
Poisson
 brackets  arising in the framework of the associativity equations
\cite{D93}, \cite{D94}. Compatible Poisson brackets of
hydrodynamic type can also be obtained as a result of  Whitham
averaging (dispersionless limit) from the local
compatible Poisson brackets of integrable systems
\cite{DN83}, \cite{DN84}, \cite{DN89}, \cite{Tsarev90}, \cite{Pavlov94},
\cite{FerPav},  \cite{Strachan}.
Some further examples and partial classification results  can be found in
\cite{Fer90}, \cite{FerPav}, \cite{Nutku}, \cite{Pavlov93},
\cite{Mokhov99}, \cite{Fordy}.

If the spectrum of $r^i_j$ is simple, the vanishing of the Nijenhuis tensor
 implies the existence of a coordinate system where both metrics $g^{ij}$
and 
$\tilde g^{ij}$ become diagonal. In these diagonal coordinates the
compatibility conditions
take the form of an integrable reduction of the Lam\'e equations. We present
the 
corresponding Lax pairs in sect. 3. Another approach to the
integrability of this system has been proposed recently by
Mokhov \cite{Mokhov001} by an appropriate modification of Zakharov's scheme
\cite{Zakharov}. 

The main observation of this paper is the relationship between
compatible Poisson brackets of hydrodynamic type and
hypersurfaces $M^{n-1}\in E^n$ which possess nontrivial deformations
preserving
the Weingarten operator. For surfaces $M^2\in E^3$ these
deformations have been investigated by Finikov and
Gambier as far back as in 1933 \cite{Finikoff}, \cite{Finikoff1}.
In sect. 4 we demonstrate that the n-orthogonal coordinate system in $E^n$
corresponding to 
the flat metric $\tilde g^{ij}+\lambda g^{ij}$
(rewritten in the diagonal coordinates)
deforms with respect to $\lambda$ in such a way that the Weingarten
operators of the 
coordinate hypersurfaces are preserved up to  constant scaling factors.
In sect. 5 we  discuss surfaces $M^2\in E^3$ which possess nontrivial
one-parameter 
 deformations preserving the Weingarten operator and
explicitly introduce a spectral parameter in the corresponding Gauss-Codazzi
equations.

\section{Differential-geometric criterion of the compatibility}

 To formulate the necessary and sufficient conditions of the
compatibility  we introduce
the operator $r^i_j=\tilde g^{is}g_{sj}$, which is automatically symmetric
\begin{equation}
r^i_sg^{sj}=r^j_sg^{si},
\label{symm}
\end{equation}
so that $\tilde g^{ij}=r^i_sg^{sj}=r^j_sg^{si}=r^{ij}$. In what follows
we use the first metric $g^{ij}$ for raising and lowering the indices.

\begin{theorem}

 Hamiltonian operators (\ref{H1}), (\ref{H2}) are compatible if and
only if the following conditions are satisfied:

\medskip

1. The Nijenhuis tensor of $r^i_j$ vanishes:
\begin{equation}
N^i_{jk}=r^s_j\partial_sr^i_k  -  r^s_k\partial_sr^i_j  -
r^i_s(\partial_j r^s_k - \partial_kr^s_j)=0.
\label{Nijenhuis}
\end{equation}

\medskip

2. The metric coefficients $\tilde g^{ij}= r^{ij}$ satisfy the equations
\begin{equation}
\nabla^i\nabla^j r^{kl}+\nabla^k\nabla^l r^{ij}=
\nabla^i\nabla^k r^{jl}+\nabla^j\nabla^l r^{ik}.
\label{nabla}
\end{equation}
Here $\nabla^i=g^{is}\nabla_s$ is the covariant differentiation
corresponding to the metric
$g^{ij}$. The vanishing of the Nijenhuis tensor implies the following
expression for the
coefficients $\tilde b^{ij}_k$ in terms of $r^i_j$:
\begin{equation}
2 \tilde b^{ij}_k= \nabla^ir^j_k-\nabla^jr^i_k+\nabla_k
r^{ij}+2b^{sj}_kr^i_s
\label{F}
\end{equation}

\end{theorem}

In a somewhat different form the necessary and sufficient conditions of the
compatibility were formulated in \cite{D93}, \cite{D94}, \cite{Mokhov99},
\cite{Mokhov00}.

{\bf Remark.} The criterion of the compatibility of the Hamiltonian
operators of
hydrodynamic type resembles that of the finite-dimensional Poisson
bivectors:
two skew-symmetric Poisson bivectors $\omega^{ij}$ and $\tilde \omega^{ij}$
are compatible if and only if the Nijenhuis tensor of the corresponding
recursion operator $r^i_j=\tilde \omega^{is}\omega_{sj}$ vanishes. We
emphasize that
in our situation operator $r^i_j$ does not coincide with the recursion
operator.

\bigskip

{\centerline {Proof of Theorem 1:}}

\medskip

Recall that in terms of  $g^{ij}$ and $b^{ij}_k$ the conditions
for the operator $A$ to be Hamiltonian take the form
\begin{equation}
2b^{ki}_sg^{sj}=g^{js}\partial_sg^{ik}+g^{ks}\partial_sg^{ij}-g^{is}\partial
_sg^{kj}
\label{J1}
\end{equation}
and
\begin{equation}
g^{js}\partial_sb^{ik}_n-g^{is}\partial_sb^{jk}_n+(b^{ij}_s-b^{ji}_s)b^{sk}_
n
+b^{ik}_sb^{js}_n-b^{jk}_sb^{is}_n=0,
\label{J2}
\end{equation}
respectively (the last condition follows from the identity
$g^{il}g^{js}R^k_{nls}=0$ after rewriting it in terms of $b^{ij}_k$).
 Note that  (\ref{J1}) is equivalent to a pair of simpler conditions
$$
b^{ij}_k+b^{ji}_k=\partial_kg^{ij}, ~~~~
b^{ik}_sg^{sj}=b^{jk}_sg^{si}.
$$
To write down the compatibility conditions of (\ref{H1}) and (\ref{H2}),
 we replace $g^{ij}$ and $b^{ij}_k$ by the linear combinations
$$
g^{ij}\to \lambda g^{ij}+ {\tilde g}^{ij}, ~~~~
b^{ij}_k\to \lambda b^{ij}_k + \tilde b^{ij}_k,
$$
substitute them into (\ref{J1}), (\ref{J2}), collect the terms with
$\lambda$
(terms with $\lambda^2$ and $\lambda^0$ vanish since (\ref{H1}) and
(\ref{H2}) are
Hamiltonian), and equate them to zero. Thus, (\ref{J1}) produces the first
compatibility condition
\begin{equation}
\begin{array}{c}
2\tilde b^{ki}_sg^{sj}+2 b^{ki}_s\tilde g^{sj}= \\
\ \\
\tilde g^{js}\partial_sg^{ik}+  g^{js}\partial_s \tilde g^{ik}  +
\tilde g^{ks}\partial_sg^{ij} +  g^{ks}\partial_s \tilde g^{ij}  -
\tilde g^{is}\partial_sg^{kj}-g^{is}\partial_s \tilde g^{kj}.
\end{array}
\label{C1}
\end{equation}
Similarly, (\ref{J2}) produces the second compatibility condition
\begin{equation}
\begin{array}{c}
\tilde g^{js}\partial_sb^{ik}_n +g^{js}\partial_s \tilde b^{ik}_n  -
\tilde g^{is}\partial_sb^{jk}_n - g^{is}\partial_s \tilde b^{jk}_n  + \\
\ \\
(\tilde b^{ij}_s-\tilde b^{ji}_s)b^{sk}_n +( b^{ij}_s- b^{ji}_s)\tilde
b^{sk}_n
+\tilde b^{ik}_sb^{js}_n+ b^{ik}_s \tilde b^{js}_n  - \tilde
b^{jk}_sb^{is}_n -
b^{jk}_s\tilde b^{is}_n=0.
\end{array}
\label{C2}
\end{equation}
To simplify further calculations, it is convenient to work in the
coordinates where
the flat metric $g$ assumes the constant coefficient form $g^{ij}=const$,
so that $b^{ij}_k\equiv 0$. In these coordinates the compatibility
conditions 
(\ref{C1}), (\ref{C2}) reduce to
\begin{equation}
2\tilde b^{ki}_sg^{sj}=
g^{js}\partial_s \tilde g^{ik}  +
  g^{ks}\partial_s \tilde g^{ij}  - g^{is}\partial_s \tilde g^{kj}
\label{C3}
\end{equation}
and
\begin{equation}
g^{js}\partial_s \tilde b^{ik}_n  - g^{is}\partial_s \tilde b^{jk}_n  =0,
\label{C4}
\end{equation}
respectively.
Rewriting the left-hand side of (\ref{C3}) in the form
$
2\tilde b^{ki}_sg^{sj}=
2\tilde b^{ki}_s \tilde g^{sl}(r^{-1})^j_l
$
and substituting the expressions for $\tilde b^{ki}_s \tilde g^{sl}$
from (\ref{J1}), we arrive at
$$
(\underline{\tilde g^{ls}\partial_s \tilde g^{ik}}+
\tilde g^{ks}\partial_s \tilde g^{il}-
\tilde g^{is}\partial_s \tilde g^{kl}) (r^{-1})^j_l =
\underline{g^{js}\partial_s\tilde g^{ik}}+
g^{ks}\partial_s\tilde g^{ij}-
g^{is}\partial_s\tilde g^{kj}.
$$
Cancelling the underlined terms and substituting
$\tilde g^{ij}=r^i_sg^{sj}=r^j_sg^{si}$,
we obtain
$$
(r^s_pg^{pk}\partial_s(r^l_ng^{ni})-r^s_pg^{pi}\partial_s(r^l_ng^{nk}))
(r^{-1})^j_l=
g^{ks}\partial_s(r^j_ng^{ni})-g^{is}\partial_s(r^j_ng^{nk}),
$$
or 
$$
(r^s_pg^{pk}g^{ni}\partial_sr^l_n-r^s_pg^{pi}g^{nk}\partial_sr^l_n)
(r^{-1})^j_l= 
g^{ks}g^{ni}\partial_sr^j_n-
g^{is}g^{nk}\partial_sr^j_n.
$$
Contraction with $r^m_j$ results in
$$
g^{pk}g^{ni}r^s_p\partial_sr^m_n-g^{pi}g^{nk}r^s_p\partial_sr^m_n=
g^{ks}g^{ni}r^m_j\partial_sr^j_n-
g^{is}g^{nk}r^m_j\partial_sr^j_n
$$
which is equivalent to
$$
g^{pk}g^{ni}(
r^s_p\partial_sr^m_n-r^s_n\partial_sr^m_p-r^m_j\partial_pr^j_n+r^m_j\partial
_nr^j_p)=0,
$$
implying the vanishing of the Nijenhuis tensor.

To establish the second identity (\ref{J2}), we will make use of the
formula (\ref{F}) for the
coefficients $\tilde b^{ij}_k$ in terms of $r$,
the proof of which is included in the Appendix (note that this formula is
true in arbitrary coordinate system).
 In the coordinates where $g^{ij}=const$ we have
$g^{is}\partial_s=\nabla^i,  ~~
2 \tilde b^{ij}_k= \nabla^ir^j_k-\nabla^jr^i_k+\nabla_k r^{ij},$ so that
(\ref{C4})
takes the form
$$
\nabla^j(\underline{\nabla^ir^k_n}-\nabla^kr^i_n+\nabla_n r^{ik})-
\nabla^i(\underline{\nabla^jr^k_n}-\nabla^kr^j_n+\nabla_n r^{jk})=0.
$$
Cancellation of the underlined terms and contraction with $g^{sn}$
produces (\ref{nabla}).
This completes the proof of the Theorem.

{\bf Remark.} If the spectrum of $r^i_j$ is simple, condition (\ref{nabla})
is redundant: 
it is automatically satisfied by virtue of (\ref{Nijenhuis}) and the
flatness of both metrics
$g$ and $\tilde g$. This was the motivation for me to drop condition
(\ref{nabla})
in the compatibility criterion formulated in \cite{Fer95}. However, in this
 general form the criterion
proved to be incorrect:  recently it was pointed out by Mokhov
\cite{Mokhov00} 
that in the case when the spectrum of $r^i_j$ is not simple the vanishing of
the Nijenhuis tensor
is  no longer sufficient for the compatibility.

\section{Compatibility conditions in the diagonal form: the Lax pairs}

If the spectrum of $r^i_j$ is simple, the vanishing of the
Nijenhuis tensor implies the existence of the coordinates $R^1,\ldots, R^n$
in which the  objects $r^i_j, ~ g^{ij}, ~ \tilde g^{ij}$ become diagonal.
Moreover, the $i$-th eigenvalue of $r^i_j$ depends only on the coordinate
$R^i$, so that $r^i_j = diag(\eta_i), ~ g^{ij}= diag(g^{ii}), ~ \tilde
g^{ij} =
diag(g^{ii}\eta_i)$ where $\eta_i$ is a function of $R^i$. This is a
generalization of 
the analogous observation by Dubrovin
\cite{D93} in the particular case of compatible Poisson brackets originating
from the 
theory of the associativity equations. Introducing the Lame coefficients
$H_i$ and the rotation coefficients
$\beta_{ij}$ by the formulae
\begin{equation}
H_i=\sqrt{g_{ii}}=1/\sqrt{g^{ii}}, ~~~~ \partial_i H_j=\beta_{ij} H_i,
\label{H}
\end{equation}
we can rewrite the zero curvature conditions for the metric
$g$ in the form
\begin{equation}
\partial_k\beta_{ij}=\beta_{ik}\beta_{kj},
\label{F1}
\end{equation}
\begin{equation}
\partial_i\beta_{ij}+\partial_j\beta_{ji}+\sum_{k\ne i,
j}\beta_{ki}\beta_{kj}=0.
\label{F2}
\end{equation}
The zero curvature condition for the metric $\tilde g$
imposes the additional constraint
\begin{equation}
\eta_i\partial_i\beta_{ij}+\eta_j\partial_j\beta_{ji}+\frac{1}{2}\eta_i{'}\b
eta_{ij}
+\frac{1}{2}\eta_j{'}\beta_{ji}+\sum_{k\ne i,
j}\eta_k\beta_{ki}\beta_{kj}=0,
\label{F3}
\end{equation}
resulting from (\ref{F2}) after the  substitution of
 the rotation coefficients $\tilde
\beta_{ij}=\beta_{ij}\sqrt{\eta_i/\eta_j}$ of the metric
$\tilde g$. As can be readily seen, equations (\ref{F2}) and (\ref{F3})
already imply the compatibility,
so that in the diagonalisable case condition (\ref{nabla}) of Theorem 1 is
indeed superfluous.
Solving equations (\ref{F2}), (\ref{F3}) for $\partial_i\beta_{ij}$,
 we can rewrite (\ref{F1})--(\ref{F3}) in the form
\begin{equation}
\begin{array}{c}
\partial_k\beta_{ij}=\beta_{ik}\beta_{kj}, \\
\ \\
\partial_i\beta_{ij}=\frac{1}{2} \frac{\displaystyle
\eta_i{'}}{\displaystyle \eta_j-\eta_i}\beta_{ij}
+\frac{1}{2} \frac{\displaystyle \eta_j{'}}{\displaystyle
\eta_j-\eta_i}\beta_{ji}+\sum_{k\ne i, j}
\frac{\displaystyle \eta_k-\eta_j}{\displaystyle
\eta_j-\eta_i}\beta_{ki}\beta_{kj}.
\end{array}
\label{S}
\end{equation}
It can be  verified by a straightforward calculation that system (\ref{S})
is compatible for any choice of the functions
$\eta_i(R^i)$,
and its general solution depends on $n(n-1)$ arbitrary functions of one
variable
(indeed, one can arbitrarily prescribe the value of $\beta_{ij}$ on the j-th
coordinate line). 
Under the additional "Egorov" assumption $\beta_{ij}=\beta_{ji}$, system
(\ref{S})
reduces to the one studied by Dubrovin in \cite{D94}.
For $n\geq 3$ system (\ref{S}) is essentially nonlinear. Its integrability
follows from 
the Lax pair
\begin{equation}
\partial_j\psi_i=\beta_{ij}\psi_j, ~~~~
\partial_i\psi_i=-\frac{\eta_i{'}}{2(\lambda+\eta_i)}\psi_i-\sum_{k\ne
i}\frac{\lambda+\eta_k}
{\lambda+\eta_i}\beta_{ki}\psi_k
\label{L1}
\end{equation}
with  a spectral parameter $\lambda$ (another demonstration of the
integrability of system (\ref{L1})
has been proposed recently in \cite{Mokhov001} by an appropriate
modification of 
Zakharov's approach \cite{Zakharov}).

{\bf Remark.} In fact, the Lax pair (\ref{L1}) is gauge-equivalent to the
equations
$$
(\tilde g^{ik}+\lambda g^{ik})\partial_k\partial_j \psi+
(\tilde b^{ik}_j+\lambda b^{ik}_j)\partial_k\psi=0
$$
for the Casimirs $\int \psi dx$ of the Hamiltonian operator
 $\tilde A^{ij}+\lambda A^{ij}$.

After the gauge transformation $\psi_i= \varphi_i/\sqrt{\lambda+\eta_i}$
the Lax pair (\ref{L1}) assumes the manifestly skew-symmetric form
\begin{equation}
\partial_j\varphi_i=
\sqrt{\frac{\lambda+\eta_i}{\lambda+\eta_j}}\beta_{ij} \varphi_j, ~~~~
\partial_i\varphi_i=-\sum_{k\ne i}\sqrt{\frac{\lambda+\eta_k}
{\lambda+\eta_i}}\beta_{ki} \varphi_k,
\label{L2}
\end{equation}
which is of the type discussed in \cite{Mikhailov}.
Thus, we can introduce an orthonormal frame $\vec \varphi_1, ...,
\vec \varphi_n$ 
in the Euclidean space $E^n$ satisfying the equations
\begin{equation}
\partial_j\vec \varphi_i=
\sqrt{\frac{\lambda+\eta_i}{\lambda+\eta_j}}\beta_{ij}\vec \varphi_j, ~~~~
\partial_i\vec \varphi_i=-\sum_{k\ne i}\sqrt{\frac{\lambda+\eta_k}
{\lambda+\eta_i}}\beta_{ki}\vec \varphi_k,
~~~~
(\vec \varphi_i, \vec \varphi_j)=\delta_{ij}.
\label{frame}
\end{equation}
Let us introduce a vector $\vec r$ such that
$$
\partial_i \vec r=\frac{H_i}{\sqrt{\lambda+\eta_i}} \vec \varphi_i
$$
(the compatibility of these equations can be readily verified).
In view of the formula
$$
(\partial_i \vec r, \partial_j \vec
r)=\frac{H_i^2}{\lambda+\eta_i}\delta_{ij},
$$
the radius-vector $\vec r$ is descriptive of an n-orthogonal coordinate
system in $E^n$
corresponding to the flat metric
$$
\sum_i \frac{H_i^2}{\lambda+\eta_i}(dR^i)^2.
$$
Geometrically, $\vec \varphi_i$ are the unit vectors along the coordinate
lines of 
this n-orthogonal system.

Let us discuss in some more detail the
 case $\eta_i=const=c_i$, in which system (\ref{S}) takes the form
\begin{equation}
\begin{array}{c}
\partial_k\beta_{ij}=\beta_{ik}\beta_{kj}, \\
\ \\
\partial_i\beta_{ij}=\sum_{k\ne i, j}
\frac{\displaystyle c_k-c_j}{\displaystyle c_j-c_i}\beta_{ki}\beta_{kj}.
\end{array}
\label{S1}
\end{equation}
One can readily verify that the quantity
$$
P_i=\sum_{k\ne i}(c_k-c_i)\beta_{ki}^2
$$
is an integral of system (\ref{S1}), namely,
$\partial_jP_i=0$ for  any $i\ne j$,
so that $P_i$ is a function of $R^i$. Utilising the obvious symmetry
$
R^i\to  s_i(R^i), ~~ \beta_{ki}\to \beta_{ki}/s_i'(R^i)
$
of system (\ref{S1}), we can reduce $P_i$ to $\pm 1$ (if nonzero). Let us
consider the simplest nontrivial case
$n=3, \ P_1=P_2=1, \ P_3=-1$:
$$
\begin{array}{c}
P_1=(c_2-c_1)\beta_{21}^2+(c_3-c_1)\beta_{31}^2=1, \\
\ \\
P_2=(c_1-c_2)\beta_{12}^2+(c_3-c_2)\beta_{32}^2=1, \\
\ \\
P_3=(c_1-c_3)\beta_{13}^2+(c_2-c_3)\beta_{23}^2=-1.
\end{array}
$$
Assuming $c_3>c_2>c_1$ and introducing the parametrization
$$
\begin{array}{c}
\beta_{21}={\sin p}/{\sqrt {c_2-c_1}}, ~~~~ \beta_{31}={\cos p}/{\sqrt
{c_3-c_1}}, \\
\ \\
\beta_{12}={\sinh q}/{\sqrt {c_2-c_1}}, ~~~~ \beta_{32}={\cosh q}/{\sqrt
{c_3-c_2}}, \\
\ \\
\beta_{13}={\sin r}/{\sqrt {c_3-c_1}}, ~~~~ \beta_{23}={\cos r}/{\sqrt
{c_3-c_2}},
\end{array}
$$
we readily rewrite (\ref{S1}) in the form
$$
\begin{array}{c}
\partial_1q=\mu_1 \cos p, ~~~~ \partial_1r=-\mu_1 \sin p, \\
\ \\
\partial_2p=-\mu_2 \cosh q, ~~~~ \partial_2r=\mu_2 \sinh q, \\
\ \\
\partial_3p=\mu_3 \cos r, ~~~~ \partial_3q=\mu_3 \sin r,
\end{array}
$$
where
$$
\mu_1=\sqrt{\frac{c_3-c_2}{(c_2-c_1)(c_3-c_1)}}, ~~~
\mu_2=\sqrt{\frac{c_3-c_1}{(c_2-c_1)(c_3-c_2)}}, ~~~
\mu_3=\sqrt{\frac{c_2-c_1}{(c_3-c_1)(c_3-c_2)}}.
$$
After rescaling, this system simplifies to
\begin{equation}
\begin{array}{c}
\partial_1q= \cos p, ~~~~ \partial_1r=- \sin p, \\
\ \\
\partial_2p=- \cosh q, ~~~~ \partial_2r= \sinh q, \\
\ \\
\partial_3p= \cos r, ~~~~ \partial_3q= \sin r.
\end{array}
\label{S2}
\end{equation}
Expressing $p$ and $r$ in the form
$
p=\arccos {\partial_1q}, ~ r=\arcsin {\partial_3q},
$
we can rewrite (\ref{S2}) as a triple of pairwise commuting Monge-Ampere
equations
$$
\begin{array}{c}
\partial_1 \partial_2q=\cosh q\sqrt{1-\partial_1q^2}, \\
\ \\
\partial_1 \partial_3q=-\sqrt{1-\partial_1q^2} \sqrt{1-\partial_3q^2}, \\
\ \\
\partial_2 \partial_3q=\sinh q\sqrt{1-\partial_3q^2}.
\end{array}
$$
Similar triples of Monge-Ampere equations  were obtained in \cite{Fer97} in
the classification of quadruples
of $3\times 3$ hydrodynamic type systems which are closed under the Laplace
transformations.
However, at the moment there is no explanation of this  coincidence.

\section{Deformations of n-orthogonal coordinate systems inducing rescalings
of the
 Weingarten operators of the coordinate hypersurfaces}

We have demonstrated in sect. 3 that the radius-vector $\vec r(R^1,...,
R^n)$
of the n-orthogonal coordinate system in $E^n$ corresponding to the flat
diagonal metric
$ \sum_i \frac{H_i^2}{\lambda+\eta_i}(dR^i)^2$ satisfies the equations
$$
\partial_i \vec r=\frac{H_i}{\sqrt{\lambda+\eta_i}} \vec \varphi_i,
$$
where the infinitesimal displacements of the orthonormal frame $\vec
\varphi_i$
are governed by 
$$
\partial_j\vec \varphi_i=
\sqrt{\frac{\displaystyle \lambda+\eta_i}{\displaystyle
\lambda+\eta_j}}\beta_{ij}\vec \varphi_j, ~~~~
\partial_i\vec \varphi_i=-\sum_{k\ne i}\sqrt{\frac{\displaystyle
\lambda+\eta_k}
{\displaystyle \lambda+\eta_i}}\beta_{ki}\vec\varphi_k.
$$
Since our formulae depend on the spectral parameter, we may speak of the
"deformation" of the n-orthogonal
coordinate system with respect to $\lambda$. To investigate this deformation
in some more detail,
we fix a coordinate hypersurface $M^{n-1}\subset E^n$ (say, $R^n=const$).
Its radius-vector $\vec r$
and the unit normal $\vec \varphi_n$ satisfy the Weingarten equations
$$
\partial_i\vec \varphi_n=\frac{\beta_{ni}}{H_i}\sqrt{\lambda+\eta_n}\
\partial_i \vec r, 
~~~~ i=1,..., n-1,
$$
implying that
$$
k^i=\frac{\beta_{ni}}{H_i}\sqrt{\lambda+\eta_n}
$$
are the principal curvatures of $M^{n-1}$. Since $\eta_n$ is a constant on
$M^{n-1}$, our deformation
 preserves the Weingarten operator of $M^{n-1}$ up to a
constant scaling factor $\sqrt{\lambda+\eta_n}$ (we point out that the
curvature line parametrization
$R^1,..., R^{n-1}$ is preserved by a construction). Thus, compatible Poisson
brackets of hydrodynamic type give rise to deformations of
n-orthogonal systems in $E^n$ which, up to scaling factors, preserve the
Weingarten operators of
the coordinate hypersurfaces. If we follow the evolution of a particular
coordinate hypersurface $M^{n-1}$,
this scailing factor can be eliminated by a homothetic transformation of the
ambient space $E^n$,
so that we arrive at the nontrivial deformation of a hypersurface which
preserves the 
Weingarten operator. However, this scaling factor cannot be eliminated for
all coordinate hypersurfaces
simultaneously.

\section{Surfaces in $E^3$ which possess nontrivial deformations  preserving
the Weingarten operator}

Interestingly enough, the problem of the classification of surfaces $M^2\in
E^3$ which possess nontrivial deformations preserving the
Weingarten operator has been formulated by Finikov and Gambier as
far back as in 1933 \cite{Finikoff}, \cite{Finikoff1}. Among other results,
they demonstrated that the only
surfaces possessing 3-parameter families of such deformations are the
quadrics, conformal transforms of surfaces of
revolution and all other surfaces having the same spherical
image of curvature lines (if surfaces have the same spherical
image of curvature lines
or, equivalently, related by a Combescure transformation, they can be
deformed simultaneously).

In this section we discuss surfaces which possess 1-parameter families of
such deformations.
Let $M^2\in E^3$ be a surface parametrized by coordinates $R^1, R^2$ of
curvature lines. Let
\begin{equation}
G_{11}(dR^1)^2+G_{22}(dR^2)^2
\label{1}
\end{equation}
be its third fundamental form (or metric of the Gaussian image, which is
automatically of 
constant curvature 1). Let $k^1, k^2$ be the radii of principal curvature
satisfying the Peterson-Codazzi equations
\begin{equation}
\frac{\partial_2k^1}{k^2-k^1}=\partial_2\ln \sqrt{G_{11}}, ~~~~
\frac{\partial_1k^2}{k^1-k^2}=\partial_1\ln \sqrt{G_{22}}.
\label{PC}
\end{equation}
Suppose there exists a flat metric
\begin{equation}
g_{11}(dR^1)^2+g_{22}(dR^2)^2
\label{2}
\end{equation}
such that
\begin{equation}
G_{11}=g_{11}/\eta_1, ~~~~ G_{22}=g_{22}/\eta_2,
\label{*}
\end{equation}
where $\eta_1, \eta_2$ are functions of $R^1, R^2$, respectively. One can
readily verify that under these assumptions
the metric
\begin{equation}
\tilde G_{11}(dR^1)^2+\tilde G_{22}(dR^2)^2=
\frac{g_{11}}{\lambda +\eta_1}(dR^1)^2+\frac{g_{22}}{\lambda+\eta_2}(dR^2)^2
\label{3}
\end{equation}
has constant curvature $1$ for any $\lambda$. Since equations (\ref{PC})
are still true if 
we replace $G_{ii}$ by $\tilde G_{ii}$, we arrive at a 1-parameter family of
surfaces $M^2_{\lambda}$
with the third fundamental forms  (\ref{3}) (which depend on $\lambda$) and
the principal curvatures
$k^1, k^2$ (which are independent of $\lambda$). Hence, the Weingarten
operators of 
 surfaces  $M^2_{\lambda}$ coincide. The problem of the classification of
surfaces which possess
1-parameter families of deformations preserving the Weingarten operator is
thus reduced to the classification
of  metrics (\ref{3}) which have constant
 Gaussian curvatures  $1$ for any $\lambda$.
Any such   metric generates an infinite family of deformable surfaces whose
principal curvatures
$k^1, k^2$ satisfy (\ref{PC}). In terms of the Lame coefficients
$H_1=\sqrt{g_{11}}, \ H_2=\sqrt{g_{22}}$ and
the rotation coefficients $\beta_{12}=\partial_1H_2/H_1, \
\beta_{21}=\partial_2H_1/H_2$ our problem reduces
to the nonlinear system
\begin{equation}
\begin{array}{c}
\partial_1H_2=\beta_{12}H_1, ~~~ \partial_2H_1=\beta_{21}H_2, \\
\ \\
\partial_1\beta_{12}+\partial_2\beta_{21}=0, \\
\ \\
\eta_1 \partial_1\beta_{12}+\eta_2
\partial_2\beta_{21}+\frac{1}{2}\eta_1'\beta_{12}
+\frac{1}{2}\eta_2'\beta_{21}+H_1H_2=0,
\end{array}
\label{4}
\end{equation}
which possesses the Lax pair
$$
\begin{array}{c}
\partial_1\psi=
\left(\begin{array}{ccc}
0 & -\sqrt{\frac{\lambda + \eta_2}{\lambda+\eta_1}}\beta_{21} &
\frac{H_1}{\sqrt{\lambda+\eta_1}}\\
\sqrt{\frac{\lambda + \eta_2}{\lambda+\eta_1}}\beta_{21} & 0 & 0 \\
-\frac{H_1}{\sqrt{\lambda+\eta_1}} & 0 & 0
\end{array}\right)
\psi, \\
\ \\
\partial_2\psi=
\left(\begin{array}{ccc}
0 & \sqrt{\frac{\lambda + \eta_1}{\lambda+\eta_2}}\beta_{12} & 0\\
-\sqrt{\frac{\lambda + \eta_1}{\lambda+\eta_2}}\beta_{12} & 0 &
\frac{H_2}{\sqrt{\lambda+\eta_2}}  \\
0 & -\frac{H_2}{\sqrt{\lambda+\eta_2}} & 0
\end{array}\right)
\psi.
\end{array}
$$
Geometrically, this Lax pair governs infinitesimal displacements of
 the orthonormal frame of the orthogonal coordinate
system on the unit sphere $S^2$, corresponding to the metric (\ref{3}).
In $2\times 2$ matrices it takes the form
$$
\begin{array}{c}
2\sqrt{\lambda+\eta_1}\ \partial_1\psi=
\left(\begin{array}{ccc}
 i\sqrt{{\lambda + \eta_2}}\beta_{21} & {H_1}\\
-{H_1} & - i\sqrt{{\lambda + \eta_2}} \beta_{21}
\end{array}\right) \psi, \\
\ \\
2\sqrt{\lambda+\eta_2}\ \partial_2\psi=
i \left(\begin{array}{ccc}
 -\sqrt{{\lambda + \eta_1}}\beta_{12} & {H_2}\\
{H_2} & \sqrt{{\lambda + \eta_1}} \beta_{12}
\end{array}\right) \psi.
\end{array}
$$
{\bf Remark.}
In \cite{Fer90} we established a one-to-one correspondence between
 surfaces possessing nontrivial  deformations preserving the
Weingarten operator and multi-Hamiltonian systems of hydrodynamic type.
Indeed, let us introduce the Hamiltonian operator
$$
\delta^{ij}g^{ii}\frac{d}{dx}+b^{ij}_kR^k_x
$$
associated with  the flat diagonal metric (\ref{2}), and the nonlocal
Hamiltonian operator (see \cite{MokFer})
$$
\delta^{ij}G^{ii}\frac{d}{dx}+\tilde
b^{ij}_kR^k_x+R^i_x\left(\frac{d}{dx}\right)^{-1}R^j_x,
$$
associated with the diagonal metric (\ref{1}) of constant curvature 1
(these operators are compatible by virtue of
(\ref{*})). According to the results of Tsarev \cite{Tsarev85},
\cite{Tsarev90}, equations
(\ref{PC}) imply that the systems of hydrodynamic type
$$
R^1_t=k^1(R) R^1_x, ~~~ R^2_t=k^2(R) R^2_x
$$
are automatically  bi-Hamiltonian with the respect to both
 Hamiltonian structures. Characteristic velocities of these systems
are the radii of principal curvature of the corresponding surfaces.

The results of this section generalise in a straightforward way
to multidimensional hypersurfaces $M^{n-1}\in E^n$.

\section{Appendix: formula for $\tilde b^{ij}_k$ }

To verify  formula (\ref{F}), it suffices to check the identities
\begin{equation}
\tilde b^{ij}_k+\tilde b^{ji}_k=\partial_kr^{ij},
\label{I1}
\end{equation}
\begin{equation}
\tilde b^{ik}_sr^{sj}=\tilde b^{jk}_sr^{si}.
\label{I2}
\end{equation}
 Substituting 
the expression for the covariant derivative
$$
\nabla_kr^{ij}=\partial_kr^{ij}-b^{si}_kr^j_s-b^{sj}_kr^i_s,
$$
into (\ref{F}), we readily obtain
$$
2\tilde b^{ij}_k=(\nabla^ir^j_k-\nabla^jr^i_k+b^{sj}_kr^i_s-b^{si}_kr^j_s)+
\partial_k r^{ij},
$$
where the expression in brackets is skew-symmetric in $i, j$. This proves
(\ref{I1}).

To verify (\ref{I2}), we first rewrite it in the form
\begin{equation}
(\nabla^ir^k_s-\nabla^kr^i_s+\nabla_s
r^{ik}+\underline{2b^{lk}_sr^i_l})r^{sj}=
(\nabla^jr^k_s-\nabla^kr^j_s+\nabla_s
r^{jk}+\underline{2b^{lk}_sr^j_l})r^{si}.
\label{I3}
\end{equation}
Since 
$b^{lk}_sr^i_lr^{sj}=b^{lk}_sg_{lt}r^{ti}r^{sj}=-\Gamma^k_{ts}r^{ti}r^{sj}$,
the underlined terms cancel in view of the symmetry of $\Gamma$. Contracting
(\ref{I3})
with $g_{pj}g_{mk}g_{ni}$, we arrive at
$$
g_{pj}(\nabla_nr_{ms}-\nabla_mr_{ns})r^{sj}+r^s_p\nabla_sr_{mn}=
g_{ni}(\nabla_pr_{ms}-\nabla_mr_{ps})r^{si}+r^s_n\nabla_sr_{mp},
$$
which, by virtue of the identity $r_{ms}r^{sj}=r^s_mr^j_s$, transforms to
$$
g_{pl}r^l_s(\nabla_nr_m^s-\nabla_mr_n^s)+r^s_p\nabla_sr_{mn}=
g_{nl}r^l_s(\nabla_pr_m^s-\nabla_mr_p^s)+r^s_n\nabla_sr_{mp}.
$$
In view of the identity
$$
r^l_s(\nabla_nr_m^s-\nabla_mr_n^s)=r^s_n\nabla_sr^l_m-r^s_m\nabla_sr^l_n,
$$
manifesting the vanishing of the Nijenhuis tensor (we emphasize that in
(\ref{Nijenhuis})
partial derivatives can be replaced by covariant derivatives with respect to
any symmetric affine connection without changing the Nijenhuis tensor),
the last equation can be rewritten as follows:
$$
g_{pl}(r^s_n\nabla_sr_m^l-r^s_m\nabla_sr_n^l)+r^s_p\nabla_sr_{mn}=
g_{nl}(r^s_p\nabla_sr_m^l-r^s_m\nabla_sr_p^l)+r^s_n\nabla_sr_{mp},
$$
or
$$
r^s_n\nabla_sr_{pm}-r^s_m\nabla_sr_{pn}+r^s_p\nabla_sr_{mn}=
r^s_p\nabla_sr_{nm}-r^s_m\nabla_sr_{np}+r^s_n\nabla_sr_{mp},
$$
which is obviously an identity. This proves  formula (\ref{F}).

\section{Acknowledgements}

I would like to thank O.~I. Mokhov for the references \cite{Mokhov00},
\cite{Mokhov001}.

\end{document}